\documentclass[sigconf]{acmart}

\AtBeginDocument{%
  \providecommand\BibTeX{{%
    \normalfont B\kern-0.5em{\scshape i\kern-0.25em b}\kern-0.8em\TeX}}}

  \setcopyright{acmcopyright}
  \copyrightyear{2022}

\newtheorem{theo}{Theorem}
\newtheorem{remark}{Remark}
\newtheorem{lemma}{Lemma}

\newtheorem{prop}{Proposition}

\usepackage{float}

\usepackage{graphicx}
\copyrightyear{2022}
\acmYear{2022}
\setcopyright{acmcopyright}\acmConference[AISS 2022]{2022 4th International Conference on Advanced Information Science and System}{November 25--27, 2022}{Sanya, China}
\acmBooktitle{2022 4th International Conference on Advanced Information Science and System (AISS 2022), November 25--27, 2022, Sanya, China}
\acmPrice{15.00}
\acmDOI{10.1145/3573834.3574501}
\acmISBN{978-1-4503-9793-3/22/11}
\begin{document}

\title{An improved hybrid regularization approach for extreme learning machine}

\author{Liangjuan Zhou}
\affiliation{%
  \institution{School of Mathematics, Hunan University}
  \city{Changsha}
  \country{China}
}

\author{Wei Miao}
\email{miaow@hnu.edu.cn}
\affiliation{%
  \institution{School of Mathematics, Hunan University}
  \city{Changsha}
  \country{China}
}
\authornote{Both authors contributed equally to this research.}

\renewcommand{\shortauthors}{Zhou and Miao.}

\begin{abstract}
   Extreme learning machine (ELM) is a network model that arbitrarily initializes the first hidden layer and can be computed speedily. In order to improve the classification performance of ELM, a $\ell_2$ and $\ell_{0.5}$ regularization ELM model ($\ell_{2}$-$\ell_{0.5}$-ELM) is proposed in this paper. An iterative optimization algorithm of the fixed point contraction mapping is applied to solve the $\ell_{2}$-$\ell_{0.5}$-ELM model. The convergence and sparsity of the proposed method are discussed and analyzed under reasonable assumptions. The performance of the proposed $\ell_{2}$-$\ell_{0.5}$-ELM method is compared with BP, SVM, ELM, $\ell_{0.5}$-ELM, $\ell_{1}$-ELM, $\ell_{2}$-ELM and $\ell_{2}$-$\ell_{1}$ELM, the results show that the prediction accuracy, sparsity, and stability of the $\ell_{2}$-$\ell_{0.5}$-ELM are better than the other $7$ models.
\end{abstract}




\sloppy
\maketitle
\section{Introduction}
Feedforward neural networks(FNNs), as one of the most frequently used neural networks which can be defined mathematically as:
$$
G_N(x_i) = \sum_{i = 1}^{N} { \beta_i g(\left \langle \omega_i,x_i \right \rangle + b_i)},
$$
where $x_i = (x_{i1}, x_{i2}, \dots,x_{ip}) \in \mathbb{R}^p$ is the input, $b_i$ is the bias and $g$ is the activation function. $\left \langle \omega_i, x_i \right \rangle = \sum_{j = 1}^{p} \omega_{ij}x_{ij}$ is the euclidean inner product, $\omega_{i}=(\omega_{i1}, \omega_{i2}, \dots,\omega_{ ip}) \in \mathbb{R}^p$ are the weights connecting the input and the $i$-th hidden node, and $\beta_i \in \mathbb{R}$ are the weights connecting the $i$-th hidden and output node. In terms of the traditional learning algorithm of FNNs, all parameters in the network need to be adjusted based on specific tasks. A classical learning method is the backpropagation (BP) algorithm, which is mainly solved by gradient descent: 
$$
\min\limits_{\omega_i, \beta_i, b_i} \sum_{i=1}^{n}\Vert t_{i} - G_N(x_{i}) \Vert_2^2,
$$
where $(x_{i}, t_{i}) (i = 1,2,\dots,n)$ denotes the training samples. However, a randomized learner model, different to the traditional learning of FNNs, called as Extreme learning machine(ELM) and related algorithms were proposed by Huang\cite{HUANG2006489}. In the ELM model, $\omega_i$ and $b_i$ are randomly assigned without training, so only $\beta_i$ needs to be trained. Set $\mathbf{T} = [t_1, t_2, \dots, t_n]$ and 
 \begin{align}\label{eq:1}
\mathbf{H} = 
\begin{bmatrix}
g(\langle \omega_1, x_1 \rangle + b_1)  & \dots & g(\langle \omega_N, x_1 \rangle + b_N)\\
\vdots   &   \dots &     \vdots \\
g(\langle \omega_1, x_n \rangle + b_1) & \dots & g(\langle \omega_N, x_n \rangle + b_N) 
\end{bmatrix},
\end{align}
once the input weights and biases are specified randomly with uniform distribution in $[-c, c]$, the hidden output matrix remains unchanged during the training phase. Accordingly, the output weights could be written by utilizing the least squares method:
\begin{align}\label{eq:2}
\min\limits_{\beta \in \mathbb{R}^N} {\left\{ \left \| \mathbf{H} \beta - \mathbf{T} \right \|_2^2\right \} },
\end{align}
the solution to model \eqref{eq:2} could be written as $\beta={\mathbf{H}^ \dagger} \mathbf{T}$, where $\mathbf{H}^ \dagger$ is the Moore–Penrose generalized inverse of hidden output matrix $\mathbf{H}$\cite{201708}.
The theoretical basis for the general approximation capability of ELM networks has been proposed and established by Igelnik\cite{1995Stochastic} , where the range of randomly allocated input weights and biases are data related and assigned in a constructive mode. Consequently, the scope of parameters in the algorithm implementation should be carefully estimated for diverse datasets. On the other hand, considering the sparsity of the output parameter $\beta$ for many high-dimensional data, Cao et al.\cite{CAO20142457} proposed a $\ell_1$ regular ELM model based on the sparsity of the $\ell_1$ regularization term, which takes the following form: 

\begin{align}\label{eq:3}
\min\limits_{\beta \in \mathbb{R}^N} {\left \{ \dfrac{1}{2} \Vert \mathbf{H} \beta - \mathbf{T} \Vert_2^2 + \lambda \Vert \beta \Vert_1 \right \},
}
\end{align}
where $\lambda > 0$ is a regularization parameter and $\beta$ is the output coefficient calculated by iteration. This model is called the Lasso model, and has been studied by many scholars in recent years \cite{2011Regression}.  

For the model $\eqref{eq:2}$, Fan et al. \cite{9235314} added a $\ell_{0.5}$ regularization term to the ELM model, based on the solution generated by $\ell_{0.5}$ is sparser than the $\ell_{1}$ regularization term \cite{Xu6205396}, and the model is defined as follows: 
 \begin{align}\label{eq:4}
\min\limits_{\beta \in \mathbb{R}^N} {\left \{ \dfrac{1}{2} \Vert \mathbf{H} \beta - \mathbf{T} \Vert_2^2 + \lambda \Vert \beta \Vert_{0.5} \right \},
}
\end{align}
where $\lambda > 0$ is a regularization parameter, the model can be solved by the iterative semi-threshold algorithm \cite{Xu6205396}.

The other regularization model for model  $\eqref{eq:2}$ was about the $\ell_2$ regularization term ($\ell_2$-ELM) \cite{CAO2016546}:
\begin{align}\label{eq:5}
\min\limits_{\beta \in \mathbb{R}^N} {\left \{ \dfrac{1}{2} \Vert \mathbf{H} \beta - \mathbf{T} \Vert_2^2 + \mu \Vert \beta \Vert_2^2 \right \},
}
\end{align}
where $\mu$ is a regularization parameter, and when the expression $\mathbf{H}^T\mathbf{H} + \mu \mathbf{I}$ is invertible after choosing the parameter $\mu$, then the solution of the model \eqref{eq:5} can be written as $\beta = (\mathbf{H}^T\mathbf{H} + \mu \mathbf{I})^{-1} \mathbf{I})^{-1} \mathbf{H}^T \mathbf{T}$. 

Hai et al.\cite{Hailiang0A} proposed a $\ell_2$-$\ell_1$-ELM hybrid model by integrating the sparsity of the $\ell_1$ regularization term and the stability of the $\ell_{2}$ regularization term as follows: 
\begin{align}\label{eq:6}
\min\limits_{\beta \in \mathbb{R}^N} {\left \{ \dfrac{1}{2} \Vert \mathbf{H} \beta - \mathbf{T} \Vert_2^2 + \lambda(\gamma \Vert \beta \Vert_1 + \varepsilon \Vert \beta \Vert_2^2) \right \}},
\end{align}
where $ \lambda \ge0$, $ \gamma \ge 0 $ and $ \varepsilon \ge0$ are regularization parameters. Inspired by the $\ell_{2}$-$\ell_{1}$-ELM model, according to Xu et al.\cite{XU20121225}, they found that the sparsity of the solution of the $\ell_p(p \in (0,1))$ regularization term: when $0< p < 0.5$, there is no significant difference in the sparse effect of $\ell_{p}$; when $0.5 < p < 1$, the smaller $p$, the better the sparse effect, so the $\ell_{0.5}$ regularization term can be used as a representative element of $\ell_{p}(p \in(0,1))$; Therefore, we propose the $\ell_{2}$-$\ell_{0.5}$-ELM model by combining the stability of $\ell_{2}$ regularization term and the sparsity of $\ell_{0.5}$ which is sparser than $\ell_{1}$, the new model is described as:
\begin{align}\label{eq:7}
\min\limits_{\beta \in \mathbb{R}^N} {\left \{ \dfrac{1}{2} \Vert \mathbf{H} \beta - \mathbf{T} \Vert_2^2 + \lambda(\gamma \Vert \beta \Vert_{0.5} + \varepsilon\Vert \beta \Vert_2^2) \right \}},
\end{align}

where the parameters have the same meaning as the expression of \eqref{eq:6}. The thought of adding $\ell_{0.5}$ and $\ell_2$ penalties simultaneously in the optimization model could be found in classification \cite{2006Gene,2018Hybrid}. This study mainly establishes an iterative algorithm and studies some properties of randomized learner model as Hai\cite{Hailiang0A}. In particular, we integrate the features of ELM and propose an iterative strategy for solving the hybrid model \eqref{eq:7}. The main contributions of this paper can be summarized as follows:

(i) The whole model is a non-convex, non-smooth and non-Lipschitz optimization problem due to the existence of $\ell_{0.5}$ norm. We propose a new algorithm called as an $\ell_2$-$\ell_{0.5}$-ELM algorithm. This algorithm is proved to be effective by analyzing the sum minimization problem of two convex functions with certain characteristics.

(ii) The key theoretical properties such as convergence, sparsity are derived to guarantee the feasibility of the proposed method.

(iii) Numerous experiments were carried out, including some UCI datasets collected from experts and intelligent systems fields,  gene datasets and ORL face image datasets. Experimental results show that the better performance of the proposed $\ell_{2}$- $\ell_{0.5}$-ELM algorithm. 

The rest of this paper is organized as follows. Section $2$ reviews some basic concepts and theories. Section $3$ demonstrates the iterative method by a fixed point equation and proposes a algorithm for $\ell_2$ - $\ell_{0.5}$-ELM model. In Section $4$, some theoretical results about convergence and sparsity are analyzed. In Section $5$, experimental results on UCI datasets, gene datasets and ORL face image datasets are shown. The conclusion is drawn in Section $6$.

\section{Preliminaries}
In this section, we present some fundamental concepts and convex optimization theorems primarily. Initially, it is about the half-thresholding function\cite{Xu6205396}. $\mathscr{P}(\lambda, t):\mathbb{R} \rightarrow \mathbb{R}, \lambda > 0$, which can be written as:
\begin{align}
\label{eq:half1111}
\mathscr{P}(\lambda, t) = \begin{cases}
\frac{2}{3}t\left(1+ \cos\left(\frac{2(\pi- \phi(t))}{3}\right)\right) & \vert t \vert > \frac{{3}}{4} {\lambda}^{\frac{2}{3}}\\
0 & \vert t \vert \le \frac{{3}}{4} {\lambda}^{\frac{2}{3}}\\
\end{cases}, 
\end{align}
where $\phi_(t)= \arccos\left(\frac{\lambda}{8}(\frac{|t|}{3})^{-\frac{3}{2}}\right), \pi= 3.14$, and then the corresponding half-thresholding operator ${\rm half}(\lambda, \beta):\mathbb{R}^N \rightarrow \mathbb{R}^N$ acts component-wise as:
\begin{align}\label{eq:8}
\left[{\rm half}(\lambda, \beta) \right]_i = \mathscr{P}(\lambda, \beta_i). 
\end{align}
Next, we introduce one key characteristic of the half-thresholding operator \cite{Xu6205396,Combettes2005SignalRB}:
\begin{align}\label{eq:9}
\Vert {\rm half}(\lambda, t) - {\rm half}(\lambda, t') \Vert \le \Vert t - t'\Vert.
\end{align}
Another crucial notion of convex optimization is the proximity operator \cite{Micchelli}:
$$
{\rm prox}_{\varphi} \beta = \arg \min \left \{ \dfrac{}{} \left \|u - \beta \right \|_2^2 + \varphi(u) \right \},
$$
where $\phi$ is a real-valued convex function on $\mathbb{R}^N$. A primary property of the proximity operator is drawn in Proposition \ref{Proposition 1}\cite{Combettes2005SignalRB}, which will be utilized to prove our major result.
\begin{prop}\label{Proposition 1}
Let $\varphi$ be a real-valued convex function on $\mathbb{R}^N$. Suppose $\psi(\cdot) = \varphi + \frac{\rho}{2} \| \cdot \|_2^2 + \langle \cdot, u \rangle + \sigma$, where $u \in {\mathbb{R}}^{N}$, $\rho \in [0, \infty )$, $\sigma \in \mathbb{R}$, then   
\begin{align}\label{eq:10}
{\rm prox}_{\psi} \beta = {\rm prox}_{\varphi / (1 + \rho)}((\beta - u) / (1 + \rho)). 
\end{align}
\end{prop}

\section{Solution: Fixed point iterative algorithm for the model}
For the ELM, the output matrix $\mathbf{H}$ is a bounded linear operator from $\mathbb{R}^N$ to $\mathbb{R}^m$ owing to the activation function $g(\cdot)\in (0, 1)$, which is finite. In order to further improve the accuracy and sparsity, we employ the regularization model \eqref{eq:7} to estimate the output weights of the network. We define concisely as:
$$
p_{\gamma, \varepsilon} = \gamma \| \beta \|_{0.5} + \varepsilon \| \beta \|_2^2,
$$
where $\varepsilon$, $\gamma \ge 0$, $p_{\gamma, \varepsilon}: \mathbb{R}^N \rightarrow [0, \infty)$. Then the model \eqref{eq:7} can be redefined as
\begin{align}\label{eq:11}
\min_{\beta \in \mathbf{R}^N} \left \{ \dfrac{1}{2} \| \mathbf{H} \beta - \mathbf{T} \|_2^2 + \lambda p_{\gamma, \varepsilon} \right \}.
\end{align}
Furthermore, we introduce the following Lemma and Theorem which will be utilized to solve our model:

\begin{lemma}\label{lemma:1}
For all $\lambda > 0 $ and $\beta \in \mathbb{R}^N$,the half-thresholding operator $(8)$ can be described as:
$$
{\rm half}(\lambda, \beta) = \arg \min_{u} \left \{ \dfrac{1}{2} \| u - \beta\|_2^2 + \lambda \| u \|_{0.5} \right \}. 
$$
\end{lemma}

\begin{lemma}\label{lemma:2}
For all $\lambda > 0, \gamma \ge 0, \varepsilon \ge 0$ and $\beta \in \mathbb{R}^N$, ${\rm half}({\frac{\lambda \gamma}{1 + 2 \varepsilon \lambda}}, \frac{\beta}{1 + 2 \varepsilon \lambda})$ is the proximity operator of $\lambda p_{\gamma,\varepsilon}(\beta)$.
\end{lemma}

\begin{theo}\label{theo:1}
Let $\lambda > 0$, $\gamma \ge 0$, $\varepsilon \ge 0$ and $\delta \in (0, \infty)$. Then $\beta$ is a minimizer of function \eqref{eq:11} if and only if it meets the fixed point equation:
\begin{align}\label{eq:13}
\beta = {\rm half}\left ({\frac{\delta \lambda \gamma}{1+2\varepsilon \lambda \delta}}, \frac{(\mathbf{I} - \delta \mathbf{H}^T \mathbf{H})\beta - \delta \mathbf{H}^T \mathbf{T}}{1 + 2 \varepsilon \lambda \delta} \right),
\end{align}
where the unit operator $\mathbf{I}: \mathbb{R}^N \rightarrow \mathbb{R}^N$, the definition of  $\mathbf{H}$ is shown in \eqref{eq:1}, and $\mathbf{H}^T$ represents the adjoint of  $\mathbf{H}$.
\end{theo}

Moreover, from the property of the proximity operator, we can drive a precise statement for the Lipschitz constant of a contractive map and the corresponding theorem as follows. 

\begin{theo}\label{theo:2}
Set $\lambda > 0, \gamma \ge 0, \varepsilon \ge 0$ and $\delta \in (0, \infty)$. Suppose that there exist two positive constants $\kappa_0$ and $\kappa$, such that the norm of the output matrix $\mathbf{H}$ shown in \eqref{eq:1} of the hidden layer is finite by them, namely $\kappa_0 \le \| \mathbf{H}^T \mathbf{H}\|_2 \le \kappa $, Thus $\beta$ is a minimizer of \eqref{eq:11} if and only if it is a fixed point of the Lipchitz map $\Gamma: \mathbb{R}^N \rightarrow \mathbb{R}^N$, that is, $\beta = \Gamma \beta$ where
\begin{align}\label{eq:13}
\Gamma \beta = {\rm half}\left({\frac{\delta \lambda \gamma}{1 + 2\varepsilon \lambda \delta}}, \frac{(\mathbf{I} - \delta \mathbf{H}^T \mathbf{H}) \beta + \delta \mathbf{H}^T \mathbf{T}}{1 + 2 \varepsilon \lambda \delta}\right).
\end{align}
Selecting $\delta = \frac{2}{\kappa_0 + \kappa}$, the Lipschitz constant is finite by
$q = 1 - \dfrac{2 \kappa_0}{\kappa + \kappa_0} \le 1$. In particular, if $\kappa_0 > 0$, we can get $\Gamma$ is a contractive map.
\end{theo}

Theorem \ref{theo:1} and Theorem \ref{theo:2} illustrate that the problem of $\ell_2$-$\ell_{0.5}$-ELM can be described as a fixed point algorithm. Furthermore, the next theorem will introduce the iterative procedure of the $\ell_2$-$\ell_{0.5}$-ELM.

\begin{theo}\label{theo:3}
Suppose that $\kappa_0$ and $\kappa$ are positive constants, such that the norm of the output matrix $\mathbf{H}$ shown in $(1)$ of the hidden layer is finite by them, namely, $\kappa_0 \le \| \mathbf{H}^T \mathbf{H}\|_2 \le \kappa $, and the sequence $\left \{ \beta \right \}_{l = 0}^{\infty} \subseteq \mathbf{R}^N$ is described iteratively as
\begin{align}\label{eq:14}
\beta_{l} = {\rm half} \left({\frac{\delta \lambda \gamma}{1 + 2\varepsilon \lambda \delta}}, \frac{(\mathbf{I} - \delta \mathbf{H}^N \mathbf{H})\beta_{l-1} - \delta \mathbf{H}^T \mathbf{T}}{1 + 2 \varepsilon \lambda \delta } \right),
\end{align}
where $l = 1, 2, 3, \dots, \lambda > 0, \varepsilon > 0, \gamma \ge 0$ and $\delta = \frac{2}{\kappa + \kappa_0}$. Thus $\{ \beta_{l} \}_{l = 0}^{\infty}$ strongly converges the minimizer of model $(10)$ in spite of the choice of $\beta_0$.
\end{theo}

\begin{remark} \label{remark:1}
It is not difficult to obtain from the proof of Theorem \ref{theo:3}.
$$
\Vert \beta_{l} - \beta^{*} \Vert_2 \le \dfrac{\kappa + \kappa_0}{\kappa_0 (\kappa + \kappa_0 + 4 \varepsilon \lambda )}  \left ( \dfrac{\kappa - \kappa_0}{\kappa + \kappa_0} \right )^{l} \Vert \mathbf{H}^T \mathbf{T}\Vert_2.
$$
Therefore, for each $\xi > 0$, if
$$
\dfrac{\kappa + \kappa_0}{\kappa_0 (\kappa + \kappa_0 + 4 \varepsilon \lambda )}  \left ( \dfrac{\kappa - \kappa_0}{\kappa + \kappa_0} \right )^{l}\Vert \beta_{1} - \beta_{0}\Vert_2 < \xi.
$$
namely,
$$
l > \dfrac{ \log \left(\frac{ \Vert \beta_1 - \beta_0\Vert_2 (\kappa + \kappa_0)}{\xi\kappa_0 (\kappa + \kappa_0 + 4 \varepsilon \lambda )} \right)}{\log \left( \frac{\kappa + \kappa_0}{\kappa - \kappa_0} \right)},
$$
thus
$$
\Vert \beta_{l} - \beta^{*} \Vert_2 < \xi.
$$
\end{remark}

As a conclusion, the complete $\ell_2$-$\ell_{0.5}$-ELM algorithm is given in Algorithm $1$ which integrates the result of Theorem \ref{theo:3} and Remark \ref{remark:1}. Next section, we want give some properties of our proposed algorithm.
\begin{table}[h]
    \centering
    \begin{tabular}{p{19.6pc}}
        \toprule
        \textbf{Algorithm 1:} the algorithm for $\ell_2$-$\ell_{0.5}$-ELM model \\
        \midrule 
        \textbf{Input:}Given a set of training samples $\mathscr{f} = \left \{ (x_j, t_j): x_j \in \mathbb{R}^{p}, t_j \in \mathbb{R}^m, j = 1, 2, \dots, n \right \}$, activation function $g$, hidden node number $N$, the related regularization parameters $\lambda > 0$, $\gamma \ge 0$, $\varepsilon \ge 0$, the corresponding parameter $\delta$, and an acceptable error $\xi$.\\
        \textbf{Step 1:} Randomly assign  a proper scope for input weight $\omega_i$ and bias $b_i(i = 1, 2,\dots, N)$ \\
        \textbf{Step 2:} Compute the hidden layer output matrix $\mathbf{H}$;\\
        \textbf{Step 3:} Set $\beta_0 = \left(0, 0, \dots, 0 \right)$, $\beta_1 = {\rm half}({\dfrac{\delta \lambda \gamma}{1+2 \varepsilon \lambda \delta}}, \dfrac{(\mathbf{I} - \delta \mathbf{H}^T \mathbf{H}) \beta_0 + \delta \mathbf{H}^T \mathbf{T}}{1 + 2 \varepsilon \lambda \delta})$, and $l_{max}$ be a minimal positive integer  but larger than
        $\dfrac{ \log \left(\dfrac{ \Vert \beta_1 - \beta_0\Vert_2 (\kappa + \kappa_0)}{\xi\kappa_0 (\kappa + \kappa_0 + 4 \varepsilon \lambda )} \right)}{\log \left( \frac{\kappa + \kappa_0}{\kappa - \kappa_0} \right)}. $ \\
        \textbf{Step 4:} For $l = 1:l_{max}$ \\
        \quad\quad\quad\quad if $l \ge l_{max}$, stop;\\ \quad\quad\quad\quad else $l:=l+1$ and update the $\beta$ as follows:
        $\beta_{l+1} = {\rm half}({\dfrac{\delta \lambda \gamma}{1+2 \varepsilon \lambda \delta}}, \dfrac{(\mathbf{I} - \delta \mathbf{H}^T \mathbf{H}) \beta_l + \delta \mathbf{H}^T \mathbf{T}}{1 + 2 \varepsilon \lambda \delta}). $\\
        repeat \textbf{Step 4}, until that the desired output weight is $\hat{\beta} = \beta_{max}$.\\
        \textbf{Output:} Return the output weights $\hat{\beta}$;\\
        \bottomrule
    \end{tabular}%
    \label{tab:addlabel}%
\end{table}

\section{Some characteristics for \texorpdfstring{$\ell_{2}$-$\ell_{0.5}$-ELM}{PDFstring}}
For the new section, we want to discuss and analyze some key characteristics of the estimator regarding $\ell_2$-$\ell_{0.5}$-ELM, such as the convergence and sparsity. 

\begin{theo}
$\beta_{l}$ strongly converges to the minimum value $\beta^{*}$ of the minimization problem
$$
\min_{\beta \in \mathbf{R}^N} \left \{ \dfrac{1}{2} \| \mathbf{H}\beta - \mathbf{T}\|_2^2 + \lambda p_{\gamma \varepsilon}(\beta)
\right \}
$$
as $l \rightarrow \infty$.
\end{theo}

$\beta_{0.5}$ in the $\ell_2$-$\ell_{0.5}$-ELM is a highly significant part of the sparsity of
the solution. Thus, we set the Theorem \ref{theo:5} as follows.

\begin{theo}\label{theo:5}
Suppose $\lambda > 0, \gamma > 0$, then the support of ${\rm half}({\frac{\lambda \gamma}{1+2 \varepsilon \lambda }}, \frac{\beta}{1 + 2 \varepsilon \lambda})$ is finite
for any $\beta \in \mathbb{R}^N$. Particularly, $\beta^{*}$ and $\beta_{l}$ are all finitely supported.
\end{theo}

If the regularization parameters $\lambda$ and $\gamma$ are fixed as some constant values, then $\beta^{*}$ and $\beta_{l}$ have only a few finite nonzero coefficients, and hence the solution to \eqref{eq:11} is sparse.

\section{Performance evaluation}
In the new section, a succession of experiments, containing some UCI benchmark datasets\cite{Hailiang0A} and gene data, are carried out to demonstrate the performance of the proposed $\ell_2$-$\ell_{0.5}$-ELM method. All the experiments are performed in the Mac Pycharm environment running on Quad-Core Intel Core i5, CPU (8 GB 2133 MHz LPDDR3) processor with the speed of 1.40GHz. The activation function of networks used in the experiments is taken as sigmoid function $g(x) = 1 / (1 + e ^{-x})$.

The $\ell_2$-$\ell_{0.5}$-ELM model is compared with seven other models: BP, SVM, ELM, $\ell_2$-$\ell_{1}$-ELM, $\ell_2$-ELM, $\ell_1$-ELM, $\ell_{0.5}$-ELM. BP includes only one hidden layer and output layer, and all parameters are trained by back-propagation algorithm; $\ell_1$-ELM and $\ell_{0.5}$-ELM are the simplified forms of $\ell_2$-$\ell_{1}$-ELM  and $\ell_2$-$\ell_{0.5}$-ELM, respectively. The activation function is defined as: $g(x) = 1/(1 + e^{-x})$. 

In order to check the algorithm for $\ell_2$-$\ell_{0.5}$-ELM model, three real classification datasets from the UCI machine learning repository are considered. The basic information of each dataset is shown in Table \ref{table_ELM_all}. The average of $30$ experimental validations was used as the final result. For these datasets, the sample size is fixed, but each sample is randomly assigned as training or testing data. 
\subsection{Performance for UCI datasets}
To validate the performance of the proposed $\ell_2$-$\ell_{0.5}$-ELM model, three types of real classification datasets were used for the experiments, including UCI\cite{2013UCI}, gene expression, and ORL face datasets. The UCI machine learning repository (2013UCI) contains three datasets: Austrian Credit Approval(Austrian), Ionosphere, and Balance Scale(Balance). The gene expression datasets contain colon\cite{Alon6745} and DLBCL\cite{2011The}, both of which are binary datasets. Moreover, the ORL face dataset includes $400$ images divided into $40$ categories. Each category contains $10$ images with different facial details and each image size is $112 \times 92$. The detail information of all datasets are summarized in Table \ref{table_ELM_all}. In addition, these data were obtained from different application fields, and it is hoped that the $\ell_2$-$\ell_{0.5}$-ELM model can be analyzed from multiple perspectives by using these data from different backgrounds. 

\setcounter{table}{0}
\begin{table}
\centering
\caption{Details of the $6$ datasets} 
 \label{table_ELM_all}
 \begin{tabular}{llccc}  
\toprule
Dataset & Type & Sapmple& Feature& Catagory \\
\midrule   
Austrian & UCI& 690& 14& 2 \\
Ionosphere& UCI& 151& 34& 2 \\
Balance & UCL& 625& 4& 3 \\
colon& gene& 62& 2000& 2 \\
DLBCL& gene& 77& 7129& 2\\
ORL& image& 400& 10304& 40 \\
 \bottomrule
    \end{tabular}
\end{table}

\begin{table}
\centering
\caption{Performance comparison of 8 models on 3 different datasets} 
\label{table_ELM1_TAA} 
\fontsize{7}{9}\selectfont
\begin{tabular}{cllcl} 
\toprule
Datasets    & Methods                          & Times(s)   &Remaining Nodes    & Accuracy($\% \pm{\%}$)  \\
\midrule
Austrain  & BP                          &2.1751     & 600   & 72.58 $\pm$ 13.57 \\  
          & SVM                         &\textbf{0.0448}     & ---   & 79.14 $\pm$ 1.98  \\
          & ELM                         &0.0588     &600    &65.37  $\pm$ {3.08} \\
          & $\ell_{0.5}$-ELM            &5.8542     &48.5    &82.76 $\pm$ 0.00\\
          &$\ell_{1}$-ELM               &8.1648     &118.5  &81.38 $\pm$ 0.00\\
          & $\ell_{2}$-ELM              &8.2735     &600    &80.36 $\pm$ 0.00\\
          & $\ell_{2}$-$\ell_{1}$-ELM   &10.041     &492.5  &81.38 $\pm$ 0.00\\
          & $\ell_{2}$-$\ell_{0.5}$-ELM &7.5875     &118.5  &\textbf{82.76 $\pm$ 0.00}\\
Ionosphere & BP                          &2.1751     &600    &72.58 $\pm$ 13.57 \\ 
          & SVM                         &0.0108     &--     &86.51 $\pm$ 2.09  \\
          & ELM                         &\textbf{0.0003}     &600    &91.55 $\pm$ 2.78  \\
          & $\ell_{0.5}$-ELM            &0.0487     &29.5   &96.96 $\pm$ 0.00   \\
         &$\ell_{1}$-ELM               &5.4755     &115.9  &97.24 $\pm$ 1.06    \\
        & $\ell_{2}$-ELM              &0.0520     &600    &96.05 $\pm$ 1.57  \\
          & $\ell_{2}$-$\ell_{1}$-ELM   &4.4093     &437.5  &96.84 $\pm$ 0.98   \\
          & $\ell_{2}$-$\ell_{0.5}$-ELM &0.0569     &193    &\textbf{98.01 $\pm$ 0.00}  \\
Balance   &BP                             &4.3814    &600       &59.99 $\pm$ 25.26   \\            
          &SVM                            &0.0215    &--        &88.63 $\pm$ 1.86   \\   
          & EL,M                           &\textbf{0.0008}    &600       &50.72 $\pm$ 6.66   \\
          & $\ell_{0.5}$-ELM              &0.1285    &23.3      &90.55 $\pm$ 0.00   \\
          &$\ell_{1}$-ELM                 &6.5074    &42.9      &90.47 $\pm$ 1.66   \\
           & $\ell_{2}$-ELM                &0.1579    &600       &90.55 $\pm$ 0.00   \\
          & $\ell_{2}$-$\ell_{1}$-ELM     &6.8678    &246.4     &90.10 $\pm$ 1.35   \\
          & $\ell_{2}$-$\ell_{0.5}$-ELM   &0.0974    &52.7      &\textbf{90.91 $\pm$ 0.00}  \\
\bottomrule
\end{tabular}
\end{table}
We repeat $30$ trials and take the averages as the final results on account of reducing the random error. 
And the regularization parameters are used to control the trade-off between the error and the penalty. 
For Austrian dataset, take the parameters ( $\ell_2$-$\ell_{0.5}$-ELM, $\ell_2$-$\ell_1$-ELM : $\lambda = 0.8, \gamma = 0.1 , \varepsilon=0.9$) and for Ionosphere dataset, take ( $\ell_2$-$\ell_{0.5}$-ELM, $\ell_2$-$\ell_1$-ELM : $\lambda = 0.9, \gamma = 0.05, \varepsilon=0.9$) and Balance Scale dataset, ( $\ell_2$-$\ell_{0.5}$-ELM : $\lambda =0.8 , \gamma = 1, \varepsilon=1$, for $\ell_2$-$\ell_1$-ELM : $\lambda = 0.005, \gamma =0.5, \varepsilon=0.5$), we set the acceptable error $\xi = 0.0001, 0.001, 0.0001$ respectively. The number of hidden nodes in the experiments is $600$. Table \ref{table_ELM1_TAA} shows the running time, the number of nodes retained, and the accuracy of the test for each dataset for the eight models (the standard deviation is kept to $4$ significant digits, $0.00$ in the table indicates a standard deviation of less than $10^{-4}$). These indices are used to measure the sparsity, stability and effectiveness of the proposed method. The corresponding figures on testing are shown as follows.

\begin{figure*}
\centering
\includegraphics[scale=0.4]{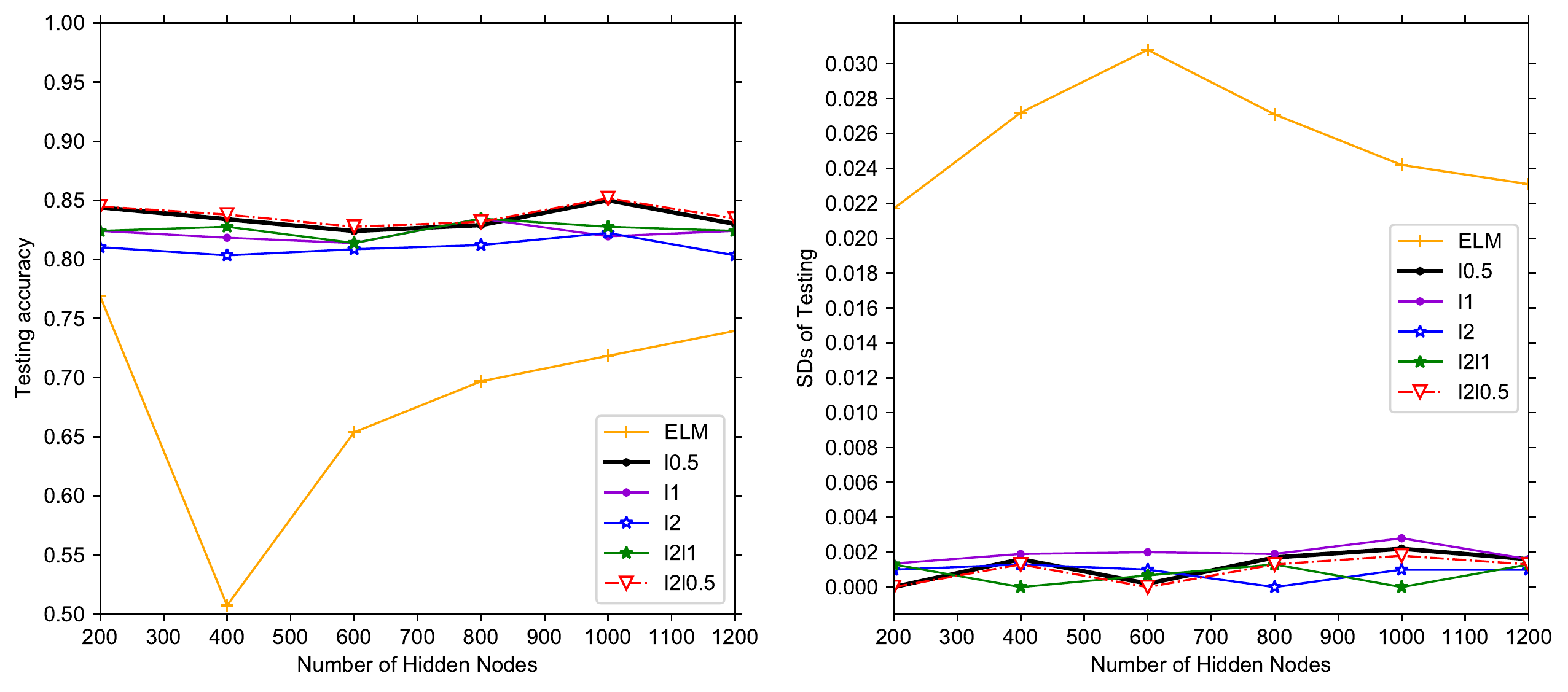}
\caption{Performance comparison of $6$ models in the Austrian dataset}
\label{fig:1}
\end{figure*}

\begin{figure*}
\centering
\includegraphics[scale=0.4]{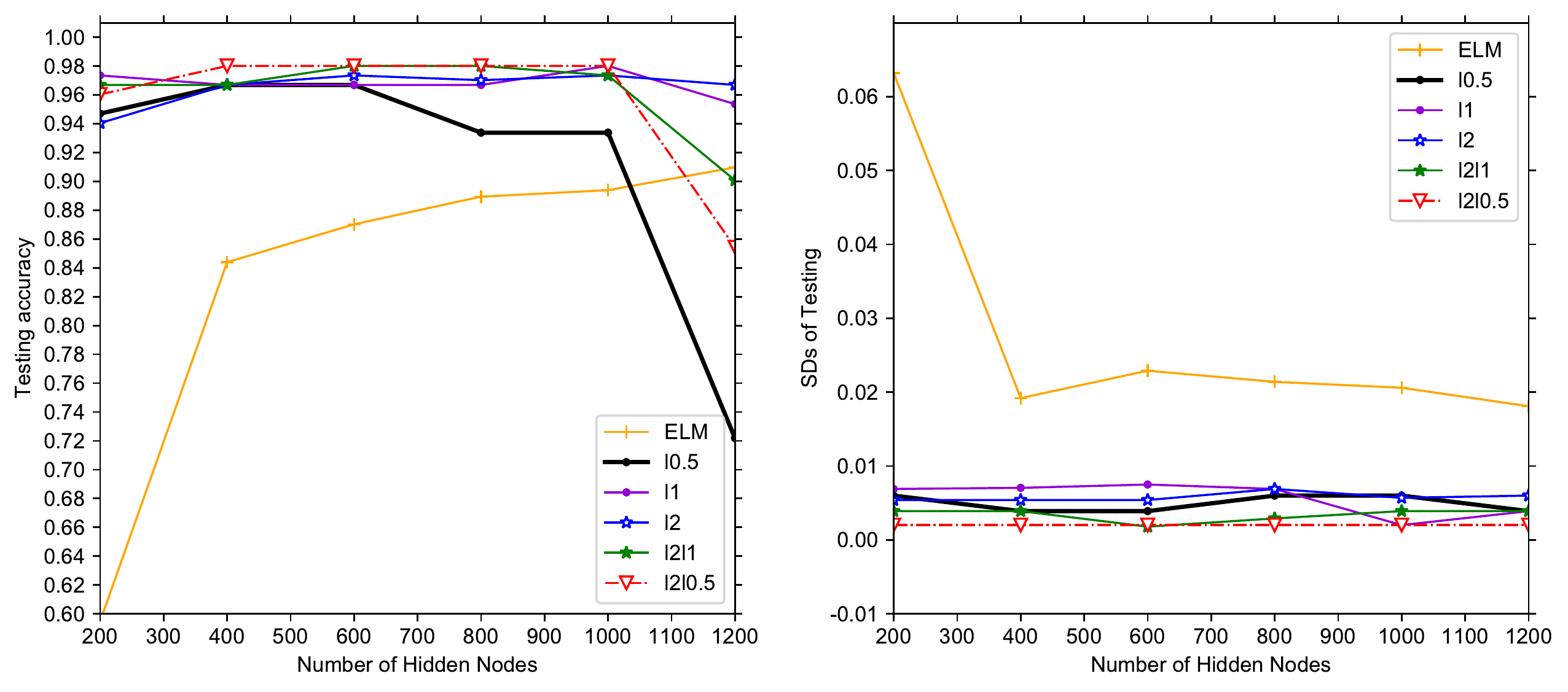}
\caption{Performance comparison of $6$ models in the Ionosphere dataset}
\label{fig:2}
\end{figure*}

\begin{figure*}
\centering
\includegraphics[scale=0.4]{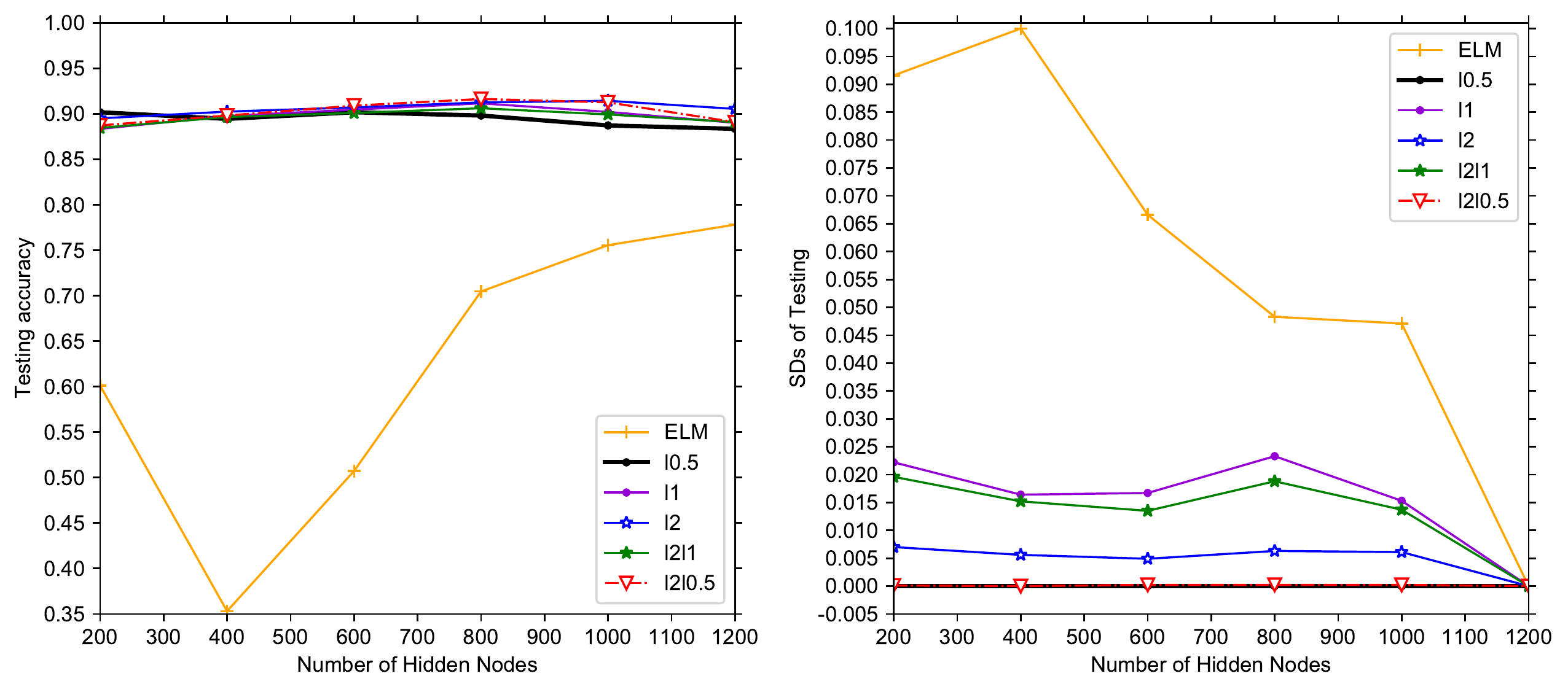}
\caption{Performance comparison of $6$ models in the Balance dataset}
\label{fig:3}
\vspace*{-8pt}\end{figure*}

From the results of \ref{fig:1}-\ref{fig:3}, we can see that the accuracy of the ELM model is lower than all the regularized ELM models. The standard deviation of the ELM model is higher than that of other regularized ELM models, which indicates that the stability of the ELM model is lower. 
The accuracy of the $\ell_2$-$\ell_{0.5}$-ELM model at all nodes can be compared with other regularized ELM models, and the accuracy at most hidden nodes is higher than other comparable regularized ELM models.  
This indicates that the $\ell_2$-$\ell_{0.5}$-ELM model has consistently good classification prediction.
In terms of the standard deviation of different nodes, the $\ell_{2}$-$\ell_{0.5}$-ELM model is lower than the other compared models, indicating that the classification accuracy of this method is more stable.

We can see the performance of $\ell_2$-$\ell_{0.5}$-ELM in detail and draw the following conclusions:

(i) In $3$ datasets, the classification accuracy of the regularized ELM methods ($\ell_2$-$\ell_{0.5}$-ELM, $\ell_{0.5}$-ELM, $\ell_2$-$\ell_{1}$-ELM, $\ell_{1}$-ELM, $\ell_{2}$-ELM) are significantly higher than that of the BP, SVM and ELM methods, indicating that the regularized ELM methods have better generalization performance, and the classification accuracy of $\ell_2$-$\ell_{0.5}$-ELM methods is higher than that of other compared regularized ELM methods.

(ii)  From the perspective of the number of remaining hidden nodes, $\ell_{0.5}$-ELM has the lowest number of hidden nodes. It is shown that the $\ell_{0.5}$ or $\ell_1$-regularization term is beneficial to enhance the sparsity of the hidden nodes of the model. Compared with the $\ell_{2}$-$\ell_{1}$-ELM model, the $\ell_{2}$-$\ell_{0.5}$-ELM model adds the $\ell_{0.5}$ regularization term to the model, which has a sparser solution and thus a better generalization ability.

(iii) From the perspective of algorithm running time, the ELM model runs in the shortest time (the ELM model can obtain the analytic solution directly without iterative computation). In comparison, the SVM model runs faster than all ELM methods with regularity. Secondly, for the $5$ regularized ELM models, the models with $\ell_{0.5}$ regularization terms ($\ell_{0.5}$-ELM, $\ell_{2}$-$\ell_{0.5}$-ELM) are faster than the models with $\ell_1$ regularization terms ($\ell_{1}$-ELM, $\ell_{2}$- $\ell_{1}$-ELM). 

\subsection{Performance for gene datasets}
In this section, the performance of the $\ell_{2}$-$\ell_{0.5}$-ELM model is validated using the colon and DLBCL data. The training and testing sets of each dataset were experimented in the ratio of $1:1$. 
The regularization parameters are set as follows, colon data: ($\ell_2$-$\ell_{0.5}$-ELM and $\ell_2$-$\ell_1$-ELM : $\lambda = 0.09, \gamma = 0.9 , \varepsilon=0.9$), DLBCL data: ($\ell_2$-$\ell_{0.5}$-ELM and $\ell_2$-$\ell_1$-ELM : $\lambda = 0.005, \gamma = 0.5, \varepsilon=0.5$); and $\xi = 0.001$. Each dataset was repeatedly run $30$ times, and the average was taken as the final result. As shown in Table \ref{table_ELM1_colon_DLBC}.

\begin{table}
\centering
\caption{Performance comparison of $8$ models in $2$ gene datasets} 
 \label{table_ELM1_colon_DLBC}
\fontsize{7.6}{9.6}\selectfont
 \begin{tabular}{clllll} 
\toprule
Datasets    & Methods  & Times(s)   &Remaining Nodes    & Accuracy($\% \pm{\%}$)  \\
\midrule
colon     &BP                          &22.2641 &1000.0  &55.52 $\pm$ 9.15  \\          
          &SVM                          &0.0358  &--  &77.5 $\pm$ 7.28   \\   
          & ELM                         &0.0056  &1000.0  &83.02 $\pm$ 1.92    \\
          & $\ell_{0.5}$-ELM            &0.0829  &370.5  &75.00 $\pm$ 0.00    \\
          &$\ell_{1}$-ELM               &0.0488  &974.5  &84.79 $\pm$ 2.22   \\
         & $\ell_{2}$-ELM               &0.0815  &1000.0  &84.17 $\pm$ 2.20  \\
          & $\ell_{2}$-$\ell_{1}$-ELM   &0.0401  &1000.0  &83.96 $\pm$ 2.24   \\
          & $\ell_{2}$-$\ell_{0.5}$-ELM &0.0879  &877.0  &\textbf{87.50  $\pm$ 0.00}   \\
DLBCL      &BP                            &122.3174   &1000.0     &57.24 $\pm$ 12.55   \\            
          &SVM                           &0.0968     &--        &87.24 $\pm$ 5.98   \\   
          & ELM                          &0.0060     &786.0     &89.90 $\pm$ 5.98   \\
          & $\ell_{0.5}$-ELM             &5.2214     &242.0     &\textbf{91.43 $\pm$ 0.00} \\ 
          &$\ell_{1}$-ELM                &18.2957    &188.5     &89.05 $\pm$ 5.12 \\
          & $\ell_{2}$-ELM               &5.2324     &764.0     &89.51 $\pm$ 5.48  \\
          & $\ell_{2}$-$\ell_{1}$-ELM    &15.5286    &431.5     &89.62 $\pm$ 6.10   \\
          & $\ell_{2}$-$\ell_{0.5}$-ELM  &5.4519     &575.5     &\textbf{91.43 $\pm$ 0.00}  \\
 \bottomrule
\end{tabular}
\end{table}
\begin{figure*}
\centering
\includegraphics[scale=0.4]{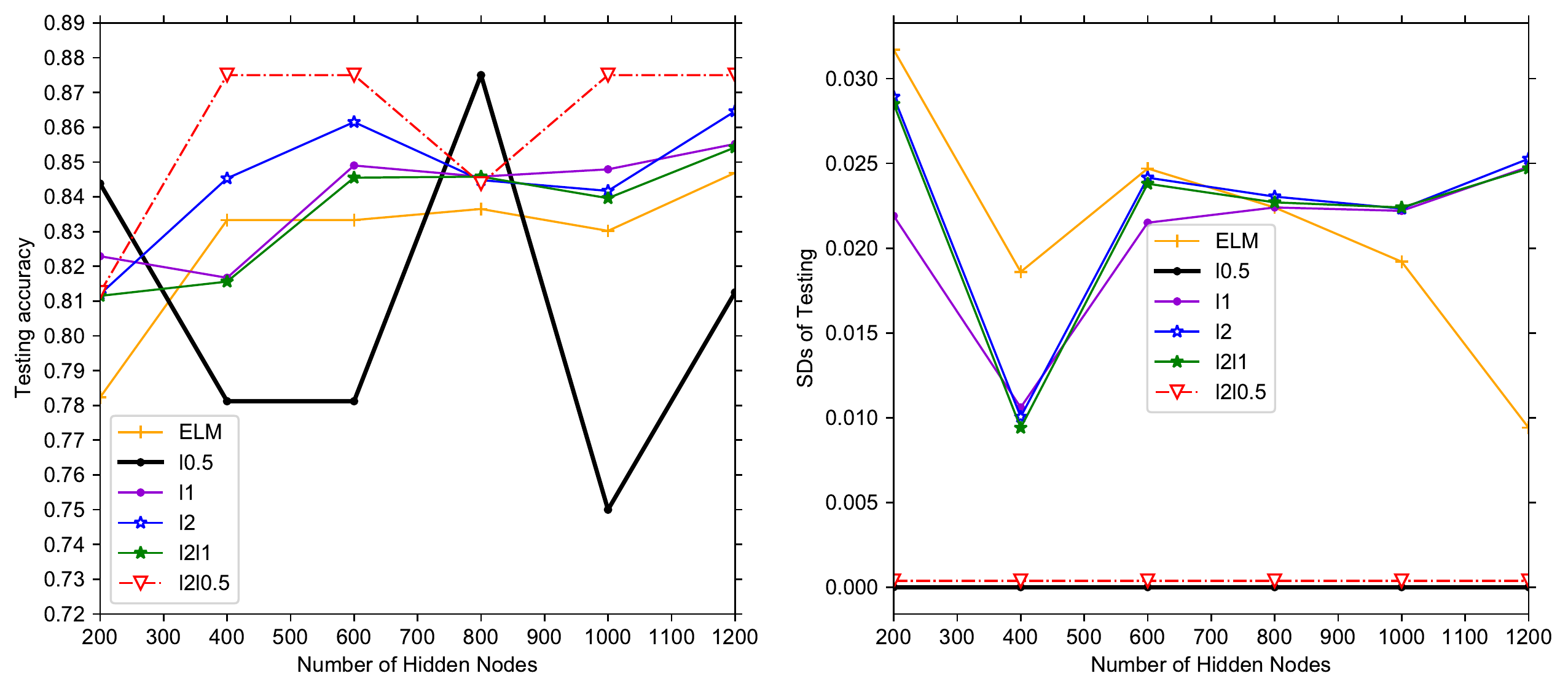}
\caption{Performance comparison of $6$ models in colon dataset}
\label{fig:4}
\end{figure*}

It can be demonstrated that the prediction accuracy of the single-layer BP network is very low and does not capture the features of the data very well. It can also be found that the prediction accuracy of the $\ell_2$-$\ell_{0.5}$-ELM model is slightly higher than that of the other methods. The standard deviations of the accuracy of the ELM methods with $\ell_{0.5}$ regularization are much smaller than those of BP, SVM, and ELM, indicating that the ELM model variants with $\ell_{0.5}$ regularization terms can improve the stability of the solutions;
 
The number of hidden nodes in the $\ell_{0.5}$-ELM and $\ell_{1}$-ELM models is smaller, that is, the sparsity of these two regularization terms is the strongest, indicating that the addition of $\ell_{0.5}$ or $\ell_1$ regularization terms in the ELM model enhances the sparsity of the model, while the number of hidden nodes in the $\ell_2$-ELM model is $1000$. The number of nodes in the $\ell_2$-ELM model is $1000$, indicating that the $\ell_2$-regularization term has no sparse effect on the model. The $\ell_2$ norm is used to increase the stability of the model by penalizing oversized regularization parameters. This makes the $\ell_{2}$-$\ell_{0.5}$-ELM sparser and model stable, and thus obtains better generalization ability.

From the perspective of algorithm running time, it can be seen that the ELM model has the shortest running time (the ELM model can obtain the analytical solution directly without iterative solving). 
In contrast, the SVM model runs faster than all ELM methods with regularization.

Further, we use the colon data to verify the effect of different number of hidden nodes ($200, 400, 600, 800, 1000, 1200$) on the stability of the ELM correlation model. We perform $30$ experiments for each hidden node and calculate the ELM, $\ell_2$-$\ell_{0.5}$-ELM, $\ell_{0.5}$-ELM, $\ell_2$-$\ell_{1}$-ELM, $\ell_{1}$-ELM, $\ell_{2}$-ELM for the test set accuracy and standard deviation as shown in Figure \ref{fig:4}. The test accuracy of $\ell_2$-$\ell_{0.5}$-ELM at all nodes can be compared with all regularized ELM models, while the accuracy at most hidden nodes is higher than other models. The standard deviation of $\ell_2$-$\ell_{0.5}$-ELM model is lower than other regularized ELM models.

\begin{table}[b]
\centering 
\caption{ Performance comparison of $8$ models in ORL face dataset} 
 \label{table_ELM_ORL_all}
 \begin{tabular}{clllll} 
\toprule
 & Methods                       & Accuracy($\%$)  \\
\midrule  
          &BP                           &31.00 $\pm$ 4.90   \\            
          &SVM                          &71.53 $\pm$ 2.12   \\   
          & ELM                         &70.58 $\pm$ 2.95   \\
          & $\ell_{0.5}$-ELM            &71.00 $\pm$ 2.34    \\
          &$\ell_{1}$-ELM               &70.85 $\pm$ 2.86   \\
          & $\ell_{2}$-ELM              &71.17 $\pm$ 2.47  \\
          & $\ell_{2}$-$\ell_{1}$-ELM   &70.58 $\pm$ 2.87   \\
          & $\ell_{2}$-$\ell_{0.5}$-ELM &\textbf{71.67 $\pm$ 2.34}\\
\bottomrule
\end{tabular}
\end{table}

\begin{table*}
    \caption{ Performance comparison of $6$ models in ORL face dataset}
    \label{table:ORL} 
\begin{tabular}{ccccccc}
    \toprule
    Nodes &ELM &$\ell_{0.5}$-ELM & $\ell_{1}$-ELM & $\ell_{2}$-ELM & $\ell_2$-$\ell_1$-ELM & $\ell_{2}$-$\ell_{0.5}$-ELM\\
    \midrule

    500  &52.92$\pm$3.04    &66.10$\pm$2.55 &60.00 $\pm$1.77  &62.63$\pm$ 2.38  &59.25 $\pm$2.32 &\textbf{65.83 $\pm $ 2.46}\\
    1500 &76.08$\pm$0.73    &77.00$\pm$0.93 &76.33 $\pm$0.67  &76.75$\pm$ 0.75  &76.33 $\pm$0.76 &\textbf{77.20 $\pm${0.93}}\\
    2000 &78.25$\pm$2.00    &78.73$\pm$2.45 &78.33 $\pm$2.08  &78.63$\pm$ 2.18  &78.33 $\pm$2.08 &\textbf{78.83 $\pm$2.45}\\
    2500 &79.58$\pm$3.49    &79.74$\pm$3.36 &79.67 $\pm$3.44  &79.21$\pm$3.29   &79.63 $\pm$3.44 &\textbf{79.76 $\pm$ 3.26}\\
    3000 &81.50$\pm$1.98    &81.55$\pm$2.69 &81.42 $\pm$2.07  &81.45$\pm$2.39   &81.42 $\pm$2.07 &\textbf{81.58 $\pm$ 2.68}\\
    3500 &81.17$\pm$1.81    &81.13$\pm$2.22 &81.17 $\pm$1.87  &81.17$\pm$1.89   &81.17 $\pm$1.87 &\textbf{81.25 $\pm$ 2.12}\\
    4000 &82.00$\pm$1.81    &82.00$\pm$1.67 &81.92 $\pm$1.74  &81.96$\pm$1.64   &81.92 $\pm$1.74 &\textbf{82.08 $\pm$ 1.65}\\
    mean &75.22$\pm$9.12    &77.16$\pm$5.33 &76.21 $\pm$7.00  &76.62$\pm$6.21   &76.08 $\pm$7.24 &\textbf{77.26 $\pm$ 5.32}\\
\bottomrule
\end{tabular}
\end{table*}

\subsection{Performance for ORL face dataset}
The ORL face dataset is used for experimental validation. The number of hidden nodes for the experiment is $1000$. The average of $30$ experiments is used as the final result. Since the original image has high dimensionality, we preprocess each image by using the $(2D)^2$PCA\cite{5946776} dimensionality reduction technique. And the training set and test set are in the ratio of $7:3$. 
The values of the regular parameters set in the experiment are as follows: $\ell_{0.5}$-ELM and $\ell_{1}$-ELM ($\gamma = 0.05, \varepsilon = 0$), $\ell_{2}$-ELM ($\gamma = 0, \varepsilon = 0.5$), $\ell_{2}$ -$\ell_{1}$-ELM, $\ell_{2}$-$\ell_{0.5}$-ELM($\gamma = 0.05, \varepsilon = 0.5$); $\lambda = 0.001$ and $\varepsilon = 0.0001$ are chosen in all experiments. 
This experiment validates the performance of the model in terms of accuracy and standard deviation. 
The results are shown in Table \ref{table_ELM_ORL_all}. 
From the table, it can be seen that the accuracy of the $\ell_{2}$-$\ell_{0.5}$-ELM model (which is slightly higher than the SVM model) is slightly higher than all other models tested.

Further, we verify the effect of different values of hidden nodes on the prediction accuracy. The number of hidden nodes chosen in the experiment is $500$, $1000$, $1500$, $2000$, $2500$, $3000$, $3500$, $4000$. The results are shown in Table \ref{table:ORL}, which show that the test accuracy of $\ell_2$-$\ell_{0.5}$-ELM model is higher than the other comparative ELM models. The test accuracy of the ELM model fluctuates the most with the changing of the number of hidden nodes, i.e., the selection of different nodes has the greatest impact on it, indicating that the ELM model is less stable in high-dimensional data.  In contrast, the standard deviations of all the regularized ELM methods ($5.33, 7.00, 6.21, 7.24, 5.32$) are lower than those of the ELM methods, indicating that the stability of the ELM model is improved by adding the regularization term. ELM methods, indicating that the stability of the proposed method is better than the other $5$ compared to ELM methods.

\section{Conclusion}
In order to further improve the stability and generalization of the ELM model, this paper proposes a $\ell_2$-$\ell_{0.5}$-ELM model by combining the $\ell_{0.5}$ and the $\ell_2$ regularization term. The iterative algorithm is applied to solve the model with a fixed points algorithm. The convergence and sparsity of this algorithm are proved. Moreover, the proposed $\ell_{2}$-$\ell_{0.5}$-ELM model is compared with BP, SVM, ELM, $\ell_{0.5}$-ELM, $\ell_{1}$-ELM, $\ell_{2}$-ELM and $\ell_2$-ELM. $\ell_2$-$\ell_{1}$-ELM models.
Experimental comparisons on several datasets (UCI dataset, gene dataset, ORL face dataset) show that the $\ell_{2}$-$\ell_{0.5}$-ELM method outperforms the other $7$ models in terms of prediction accuracy and stability on these data. Therefore, the model can be improved as follows: the information of previously computed nodes is not used in the computation of different hidden nodes, and it can be learned from the incremental learning point of view, which can reduce the computation time to a certain extent.

\balance
\bibliographystyle{ACM-Reference-Format}
\bibliography{mybib}

\end{document}